\documentclass[12pt,a4paper]{amsart} 

\usepackage{hyperref} 
\usepackage{amssymb}
\usepackage{amsrefs}
\usepackage[all]{xy}

\allowdisplaybreaks

\title[Orbit closures of directing modules]%
  {Orbit closures of directing modules \\ 
   are regular in codimension one}

\author{Grzegorz Bobi\'nski}

\address{Faculty of Mathematics and Computer Science \\ Nicolaus 
Copernicus University \\ ul.~Chopina 12/18 \\ 87-100 Toru\'n \\ 
Poland} 

\email{gregbob@mat.uni.torun.pl} 

\let\Im=\undefined
\let\mod=\undefined
\DeclareMathOperator{\GL}{GL} %
\DeclareMathOperator{\id}{id} %
\DeclareMathOperator{\Id}{Id} %
\DeclareMathOperator{\Im}{Im} %
\DeclareMathOperator{\op}{op} %
\DeclareMathOperator{\pd}{pd} %
\DeclareMathOperator{\rk}{rk} %
\DeclareMathOperator{\add}{add} %
\DeclareMathOperator{\End}{End} %
\DeclareMathOperator{\Ext}{Ext} %
\DeclareMathOperator{\Hom}{Hom} %
\DeclareMathOperator{\Ker}{Ker} %
\DeclareMathOperator{\mod}{mod} %
\DeclareMathOperator{\rad}{rad} %
\DeclareMathOperator{\gldim}{gl.dim} %
\DeclareMathOperator{\bdim}{\mathbf{dim}} %

\newcounter{claim}[section]

\newtheorem{lemma}[claim]{Lemma}
\newtheorem{corollary}[claim]{Corollary}
\newtheorem{proposition}[claim]{Proposition}

\newtheorem*{maintheorem}{Main Theorem} 

\newcommand{\bbA}{\mathbb{A}}
\newcommand{\bbB}{\mathbb{B}}
\newcommand{\bbM}{\mathbb{M}}
\newcommand{\bbN}{\mathbb{N}}
\newcommand{\bbR}{\mathbb{R}}
\newcommand{\bbV}{\mathbb{V}}
\newcommand{\bbZ}{\mathbb{Z}}

\newcommand{\bd}{\mathbf{d}}

\newcommand{\calA}{\mathcal{A}}
\newcommand{\calB}{\mathcal{B}}

\newcommand{\calE}{\mathcal{E}}
\newcommand{\calF}{\mathcal{F}}
\newcommand{\calH}{\mathcal{H}}

\newcommand{\calO}{\mathcal{O}}

\newcommand{\calT}{\mathcal{T}}
\newcommand{\calU}{\mathcal{U}}
\newcommand{\calX}{\mathcal{X}}

\newcommand{\frakR}{\mathfrak{R}}

\newcommand{\ol}{\overline}

\begin{document}

\begin{abstract} 
We show that the orbit closure of a directing module is regular in 
codimension one. In particular, this result gives information 
about a distinguished irreducible component of a module variety. 
\end{abstract} 

\maketitle


Throughout the paper $k$ is a fixed algebraically closed field. By
$\bbZ$, $\bbN$, and $\bbN_+$, we denote the sets of integers,
nonnegative integers, and positive integers, respectively. If $i,
j \in \bbZ$, then $[i, j]$ denotes the set of all $l \in \bbZ$
such that $i \leq l \leq j$.

\section*{Introduction and the main result}

Given a finite-dimensional $k$-algebra $\Lambda$ and an element
$\bd$ of the Grothendieck group $K_0 (\Lambda)$ of the category of
$\Lambda$-modules, one defines the variety $\mod_\Lambda^{\bd}
(k)$ of $\Lambda$-modules of dimension vector $\bd$. A product
$\GL_{\bd} (k)$ of general linear groups acts on
$\mod_\Lambda^{\bd} (k)$ in such a way that the $\GL_{\bd}
(k)$-orbits correspond to the isomorphism classes of
$\Lambda$-modules of dimension vector $\bd$. A study of properties
of the module varieties is an interesting and important direction
of research in the representation theory of algebras (for some
reviews of results see~\cites{Bongartz1998, Geiss1996, Kraft1982}
and for some more specific results see~\cites{BarotSchroer2001,
Bobinski2008, BobinskiSkowronski1999b}). In particular, the author
has showed~\cite{Bobinski2002} that if $\bd$ is the dimension
vector of a directing module and $\Lambda$ is tame, then
$\mod_\Lambda^{\bd} (k)$ is normal if and only if it is
irreducible. In general, if $M$ is a directing module, then the
closure $\ol{\calO (M)}$ of the $\GL_{\bd} (k)$-orbit $\calO (M)$
of $M$ is an irreducible component of $\mod_\Lambda^{\bd} (k)$.
Thus the above result naturally rises a question about the
properties of $\ol{\calO (M)}$. We also note that $\ol{\calO (M)}$
coincides with the closure of the set modules of projective
(injective) dimension at most $1$, which is known to be an
irreducible component of a module variety in many other cases
(see~\cite{BarotSchroer2001}*{Proposition~3.1} for details).

The question about properties of $\ol{\calO (M)}$ for a directing
module $M$ is a special case of another geometric problem
investigated in representation theory of finite-dimensional
algebras, namely, study of properties of orbit closures in module
varieties (see for example~\cites{BenderBongartz2003,
Bongartz1994, BobinskiZwara2002, Zwara2002b, Zwara2005}). In
particular, Zwara and the author proved~\cite{BobinskiZwara2006}
(using, among other things, the results
of~\cite{BobinskiSkowronski1999a}) that if $M$ is an
indecomposable directing module, then $\ol{\calO (M)}$ is a normal
variety. Recall that normal varieties are regular in codimension
one, i.e.\ the set of singular points is of codimension at least
two. Thus, the following main result of the paper is the first
step in order to generalize the above result about the closures of
the indecomposable directing modules over the tame algebras to
arbitrary directing modules over arbitrary algebras.

\begin{maintheorem} 
If $M$ is a directing module, then $\ol{\calO (M)}$ is regular in
codimension one.
\end{maintheorem} 

The paper is organized as follows. In Section~\ref{sect_direct} we
recall definitions of quivers and their representations. We also
describe properties of directing modules needed it the proof of
our main result. In Section~\ref{sect_ext} we discuss
interpretations of extension groups useful in geometric
investigations. Next, in Section~\ref{sect_scheme}, we define
module schemes and some schemes connected with them. Finally, in
Section~\ref{sect_proof}, we proof the main result of the paper.

The main idea of the proof is the following. We first observe that
each minimal degeneration $N$ of a directing module $M$ (i.e.\ a
module $N$ whose orbit is maximal in $\ol{\calO (M)} \setminus
\calO (M)$) over an algebra $\Lambda$ is of the form $N = U \oplus
V$ for a short exact sequence
\[
\xi : 0 \to U \to M \to V \to 0.
\]
Now we use a connection between the tangent space $T_N
\mod_\Lambda^{\bd} (k)$ to the module variety at $N$ and the first
extension group. As a consequence, it follows that $\Ext_\Lambda^2
(V, U)$ measures a difference between $\dim \ol{\calO (M)}$ and
$\dim_k T_N \mod_\Lambda^{\bd} (k)$. On the other hand, we show
for a general minimal degeneration $N$ of $M$, that if
\[
\xi_1 : 0 \to U \to W_1 \to U \to 0 \qquad \text{and} \qquad \xi_2
: 0 \to V \to W_2 \to V \to 0
\]
are short exact sequences, then $(\xi_1, \xi_2)$ corresponds to an
element of $T_N \ol{\calO (M)}$ if only only if the sequences
$\xi_1 \circ \xi$ and $\xi \circ \xi_2$ determine the same element
in $\Ext_\Lambda^2 (V, U)$. We prove that the space of such pairs
of sequences is of codimension $\dim_k \Ext_\Lambda^2 (V, U)$ in
$\Ext_\Lambda^1 (U, U) \times \Ext_\Lambda^1 (V, V)$, and this
will finish the proof.

For a basic background on the representation theory of algebras
(in particular, on tilting theory) we refer
to~\cite{AssemSimsonSkowronski2006}. A functorial approach to
schemes in algebraic geometry used in the article is explained for
example in~\cite{EisenbudHarris2000}.

An article was written while the author was staying at University
of Bielefeld as Alexander von Humboldt Foundation fellow. The
author expresses his gratitude to Professor Ringel for his
hospitality and helpful discussions.

\section{Preliminaries on quivers and their representations}
\label{sect_direct}

In this section we present basic facts about quivers and their
representations. We also collect facts about directing modules
necessary in the proof.

By a quiver $\Delta$ we mean a finite set $\Delta_0$ of vertices
and a finite set $\Delta_1$ of arrows together with two maps $s, t
: \Delta_1 \to \Delta_0$, which assign to an arrow $\alpha \in
\Delta_1$ its starting vertex $s_\alpha$ and terminating vertex
$t_\alpha$, respectively. By a path of length $n \in \bbN_+$ in
$\Delta$ we mean a sequence $\sigma = \alpha_1 \cdots \alpha_n$ of
arrows such that $s_{\alpha_i} = t_{\alpha_{i + 1}}$ for each $i
\in [1, n - 1]$. We write $s_\sigma$ and $t_\sigma$ for
$s_{\alpha_n}$ and $t_{\alpha_1}$, respectively. Additionally, for
each vertex $x$ of $\Delta$ we introduce a path $x$ of length $0$
such that $s_x = x = t_x$.

With a quiver $\Delta$ we associate its path algebra $k \Delta$,
which as a $k$-vector space has a basis formed by all paths in
$\Delta$ and whose multiplication is induced by composition of
paths. If $\rho = \lambda_1 \sigma_1 + \cdots + \lambda_n
\sigma_n$ for scalars $\lambda_1, \ldots, \lambda_n \in k$ and
paths $\sigma_1, \ldots, \sigma_n \in x (k \Delta) y$, where $x, y
\in \Delta_0$, then we put $s_{\rho} = y$ and $t_{\rho} = x$. Such
$\rho$ is called a relation in $\Delta$ if the length of
$\sigma_i$ is at least $2$ for each $i \in [1, n]$. A set $\frakR$
of relations is called minimal if for every $\rho \in \frakR$,
$\rho$ does not belong to the ideal $\langle \frakR \setminus \{
\rho \} \rangle$ of $k \Delta$ generated by $\frakR \setminus \{
\rho \}$. A pair $(\Delta, \frakR)$ consisting of a quiver
$\Delta$ and a minimal set of relations $\frakR$, such that there
exists $n \in \bbN_+$ with the property $\sigma \in \frakR$ for
each path $\sigma$ in $\Delta$ of length at least $n$, is called a
bound quiver. If $(\Delta, \frakR)$ is a bound quiver, then the
algebra $k \Delta / \langle \frakR \rangle$ is called the path
algebra of $(\Delta, \frakR)$.

Let $R$ be a commutative $k$-algebra and $(\Delta, \frakR)$ a
bound quiver. By an $R$-representation of $\Delta$ we mean a
collection $M = (M_x, M_\alpha)_{x \in \Delta_0, \, \alpha \in
\Delta_1}$ of free $R$-modules $M_x$, $x \in \Delta_0$, of finite
rank and $R$-linear maps $M_\alpha : M_{s_\alpha} \to
M_{t_\alpha}$, $\alpha \in \Delta_1$. An $R$-representation $M$ of
$\Delta$ is called an $R$-representation of $(\Delta, \frakR)$ if
$M_\rho = 0$ for all $\rho \in \frakR$, where we put $M_x =
\Id_{M_x}$ for $x \in \Delta_0$, $M_{\sigma} = M_{\alpha_1} \cdots
M_{\alpha_n}$ for a path $\sigma = \alpha_1 \cdots \alpha_n$ with
$\alpha_1, \ldots, \alpha_n \in \Delta_1$, and
\[
M_\rho = \lambda_1 M_{\sigma_1} + \cdots + \lambda_n M_{\sigma_n}
\]
for $\rho = \lambda_1 \sigma_1 + \cdots + \lambda_n \sigma_n$ with
scalars $\lambda_1, \ldots, \lambda_n \in k$ and paths $\sigma_1,
\ldots, \sigma_n \in x (k \Delta) y$, where $x, y \in \Delta_0$.
By a morphism $f : M \to N$ we mean a collection $(f_x)_{x \in
\Delta_0}$ of $R$-linear maps $f_x : M_x \to N_x$, $x \in
\Delta_0$, such that $f_{t_\alpha} M_\alpha = N_\alpha
f_{s_\alpha}$ for each $\alpha \in \Delta_1$. If $\Lambda = k
\Delta / \langle \frakR \rangle$, then the category of
$R$-representations of $(\Delta, \frakR)$ is equivalent to the
full subcategory $\mod_\Lambda (R)$ of the category of
$\Lambda$-$R$-bimodules formed by the bimodules $M$ such that $x
M$ is a free $R$-module for each $x \in \Delta_0$ (see for
example~\cite{AssemSimsonSkowronski2006}*{Theorem~III.1.6} for
this statement in the case $R = k$). We will identify such
$\Lambda$-$R$-bimodules and $R$-representations of $(\Delta,
\frakR)$. For $M, N \in \mod_\Lambda (R)$ we denote by
$\Hom_\Lambda (M, N)$ the space of homomorphisms from $M$ to $N$
in $\mod_\Lambda (R)$. Moreover, if $\dim_k R < \infty$, then $[M,
N]$ denotes $\dim_k \Hom_\Lambda (M, N)$. Similarly, for $n \in
\bbN_+$ we denote by $\Ext_\Lambda^n (M, N)$ the $n$-th extension
group in the category of $\Lambda$-$R$-bimodules and, if $\dim_k R
< \infty$, by $[M, N]^n$ the dimension of $\Ext_\Lambda^n (M, N)$
over $k$. For an $R$-representation $M$ of $\Delta$ its dimension
vector $\bdim M \in \bbN^{\Delta_0}$ is defined by $(\bdim M)_x =
\rk_R M_x$ for $x \in \Delta_0$.

In the rest of this section we will work in the category
$\mod_\Lambda = \mod_\Lambda (k)$ of $\Lambda$-modules for the
path algebra $\Lambda$ of a bound quiver $(\Delta, \frakR)$. For
subcategories $\calA$ and $\calB$ of $\mod_\Lambda$ we denote by
$\calA \vee \calB$ the additive closure of their union, i.e.\ the
full subcategory of $\mod_\Lambda$ whose objects are $M \oplus N$
with $M \in \calA$ and $N \in \calB$.

By a path in $\mod_\Lambda$ we mean a sequence $X_0 \to X_1 \to
\cdots \to X_n$ of nonzero maps between indecomposable
$\Lambda$-modules for $n \in \bbN_+$. A $\Lambda$-module $M$ is
called directing if there exists no path of the form
\[
M' \to \cdots \to \tau X \to * \to X \to \cdots \to M''
\]
for an indecomposable nonprojective $\Lambda$-module $X$ and
indecomposable direct summands $M'$ and $M''$ of $M$, where $\tau$
denotes the Auslander--Reiten translation in $\mod_\Lambda$. If
$M$ is indecomposable, then $M$ is directing if and only if there
is no path in $\mod_\Lambda$ of the form $M \to \cdots \to M$
(see~\cite{HappelRingel1993}*{Section~1, Corollary}).

A $\Lambda$-module $T$ is called tilting if $[T, T]^1 = 0$,
$\pd_\Lambda T \leq 1$, and there exists an exact sequence of the
form
\[
0 \to \Lambda \to T' \to T'' \to 0
\]
for $T', T'' \in \add T$, where for a $\Lambda$-module $M$ we
denote by $\add M$ the full subcategory of $\mod_\Lambda$
consisting of direct sums of direct summands of $M$. Any tilting
$\Lambda$-module $T$ determines the torsion theory $(\calT,
\calF)$, where
\begin{align*}
\calF = \calF (T) & = \{ N \in \mod_\Lambda \mid [T, N] = 0 \}
\\
\intertext{and} %
\calT = \calT (T) & = \{ N \in \mod_\Lambda \mid [T, N]^1 = 0 \}.
\end{align*}
An algebra $\Lambda$ is called tilted, if there exists a
hereditary algebra $\Sigma$ and a multiplicity free tilting
$\Sigma$-module $T$ such that $\Lambda \simeq \End_\Sigma
(T)^{\op}$. It is known that if $\Lambda$ is tilted, then $\gldim
\Lambda \leq 2$ and there are no oriented cycles in $\Delta$.

Let $M$ be a directing $\Lambda$-module. Then the support of $M$
is convex~\cite{Bobinski2002}*{Lemma~1.1}, where by the support of
a $\Lambda$-module $N$ we mean the full subquiver of $\Delta$
whose vertices are all $x \in \Delta_0$ such that $N_x \neq 0$.
Moreover, a full subquiver $\Delta'$ of a quiver $\Delta$ is
called convex if for each path $x_0 \to \cdots \to x_n$, $n \in
\bbN_+$, in $\Delta$ with $x_0, x_n \in \Delta_0'$, $x_i \in
\Delta_0'$ for all $i \in [1, n - 1]$. Consequently, we may assume
that considered directing modules are sincere, i.e.\ $M_x \neq 0$
for all $x \in \Delta_0$. If this is the case, then there exists a
directing tilting $\Lambda$-module $T$, such that $\add M \subset
\add T$~\cite{Bakke1988}.

Assume now that $T$ is a directing tilting $\Lambda$-module. We
have the following properties, which are either a general
statements of tilting theory or were proved by
Bakke~\cite{Bakke1988}. First of all, $\Lambda$ is a tilted
algebra. Secondly, $\mod_\Lambda = \calF \vee \calT$. Moreover,
$\pd_\Lambda N \leq 1$ for all $N \in \calF \vee \add T$,
$\id_\Lambda N \leq 1$ for all $N \in \calT$, and $[N', N''] = 0$
and $[N'', N']^1 = 0$ for all $N' \in \calT$ and $N'' \in \calF$.
Finally, modules in $\add T$ are uniquely determined by their
dimension vectors and if $[X, T] \neq 0$ and $[T, X] \neq 0$ for
an indecomposable $\Lambda$-module $X$, then $X \in \add T$.

Now assume that $\Lambda$ is tilted (more generally, $\gldim
\Lambda \leq 2$ and there are no oriented cycles in $\Delta$). We
define the Euler bilinear form $\langle -, - \rangle :
\bbZ^{\Delta_0} \times \bbZ^{\Delta_0} \to \bbZ$ by
\[
\langle \bd', \bd'' \rangle = \sum_{x \in \Delta_0} d_x' d_x'' -
\sum_{\alpha \in \Delta_1} d_{s_\alpha}' d_{t_\alpha}'' +
\sum_{\rho \in \frakR} d_{s_\rho}' d_{t_\rho}''
\]
for $\bd', \bd'' \in \bbZ^{\Delta_0}$. It is known
(see~\cite{Bongartz1983}*{2.2}), that if $M$ and $N$ are
$\Lambda$-modules, then
\[
\langle \bdim M, \bdim N \rangle = [M, N] - [M, N]^1 + [M, N]^2.
\]
We also have the Euler characterisitic $\chi : \bbZ^{\Delta_0} \to
\bbZ$ defined by $\chi (\bd) = \langle \bd, \bd \rangle$.

\section{Interpretations of extension groups}
\label{sect_ext}

Throughout this section we assume that $\Lambda$ is the path
algebra of a bound quiver $(\Delta, \frakR)$. Moreover, $R$ is a
commutative $k$-algebra. Our aim is this section is to present
interpretations of extension groups, which are useful in geometric
considerations.

We first present a construction investigated by
Bongartz~\cite{Bongartz1994}. For two collections $U = (U_x)_{x
\in \Delta_0}$ and $V = (V_x)_{x \in \Delta_0}$ of free
$R$-modules of finite rank we introduce the following notation:
\begin{align*}
\bbV^{V, U} & = \prod_{x \in \Delta_0} \Hom_R (V_x, U_x), &
\bbA^{V, U} & = \prod_{\alpha \in \Delta_1} \Hom_R (V_{s_\alpha},
U_{t_\alpha}),
\\
\intertext{and} %
\bbR^{V, U} & = \prod_{\rho \in \frakR} \Hom_R (V_{s_\rho},
U_{t_\rho}).
\end{align*}
If $M \in \mod_\Lambda (R)$, then we denote also by $M$ the
corresponding collection $(M_x)_{x \in \Delta_0}$ of free
$R$-modules.

Fix $U, V \in \mod_\Lambda (R)$. If $Z \in \bbA^{V, U}$, then we
define $Z_\rho^{V, U}$ for $\rho \in x (k \Delta) y$, where $x$
and $y$ run through $\Delta_0$, in the following way: $Z_x^{V, U}
= 0$ for $x \in \Delta_0$,
\[
Z_\sigma^{V, U} = \sum_{i \in [1, n]} U_{\alpha_1} \cdots
U_{\alpha_{i - 1}} Z_{\alpha_i} V_{\alpha_{i + 1}} \cdots
V_{\alpha_n}
\]
for a path $\sigma = \alpha_1 \cdots \alpha_n$ with $\alpha_1,
\ldots, \alpha_n \in \Delta_1$, and
\[
Z_\rho^{V, U} = \lambda_1 Z_{\sigma_1}^{V, U} + \cdots + \lambda_n
Z_{\sigma_n}^{V, U}
\]
for $\rho = \lambda_1 \sigma_1 + \cdots + \lambda_n \sigma_n$ with
scalars $\lambda_1, \ldots, \lambda_n \in k$ and paths $\sigma_1,
\ldots, \sigma_n \in x (k \Delta) y$, where $x, y \in \Delta_0$.
Observe that $Z_{\rho_1 \rho_2}^{V, U} = Z_{\rho_1}^{V, U}
V_{\rho_2} + U_{\rho_1} Z_{\rho_2}^{V, U}$ for all possible
$\rho_1$ and $\rho_2$ with $s_{\rho_1} = t_{\rho_2}$. We define
$\bbZ^{V, U}$ as the set of all $Z \in \bbA^{V, U}$ such that
$Z_\rho^{V, U} = 0$ for all $\rho \in \frakR$. For $Z \in \bbZ^{V,
U}$ we define an $R$-representation $W^Z$ of $\Delta$ by $W_x^Z =
U_x \oplus V_x$ for $x \in \Delta_0$ and
\[
W_\alpha^Z =
\begin{bmatrix}
U_\alpha & Z_\alpha \\ 0 & V_\alpha
\end{bmatrix}, \; \alpha \in \Delta_1.
\]
Then $W_\rho^Z = 0$ for all $\rho \in \frakR$, hence $W^Z \in
\mod_\Lambda (R)$. Moreover, we have a short exact sequence
\[
\xi^Z : 0 \to U \xrightarrow{f^Z} W^Z \xrightarrow{g^Z} V \to 0
\]
in $\mod_\Lambda (R)$, where the maps $f^Z$ and $g^Z$ are the
canonical injection and the canonical projection, respectively. On
the other hand, for every short exact sequence
\[
\xi : 0 \to U \to M \to V \to 0
\]
of $\Lambda$-$R$-bimodules, there exists (nonunique) $Z \in
\bbZ^{V, U}$ such that $[\xi^Z] = [\xi]$. More precisely, the map
$\bbZ^{V, U} \to \Ext_\Lambda^1 (V, U)$, $Z \mapsto [\xi^Z]$, is a
surjective $R$-linear map. The kernel $\bbB^{V, U}$ of this map
consists of $Z \in \bbZ^{V, U}$ such that the corresponding
sequence splits, i.e.\ there exists $h \in \bbV^{V, U}$ such that
$Z_\alpha = U_\alpha h_{s_\alpha} - h_{t_\alpha} V_\alpha$ for all
$\alpha \in \Delta_1$. From now on we identity $\Ext_\Lambda^1 (V,
U)$ with $\bbZ^{V, U} / \bbB^{V, U}$. Observe that
\[
\dim_k \bbZ^{V, U} = [V, U]^1 - [V, U] + \dim_k R \cdot \sum_{x
\in \Delta_0} \rk_R U_x \cdot \rk_R V_x
\]
provided $\dim_k R < \infty$.

We also present interpretations of some other homological
constructions.

\begin{proposition} \label{prop_exthom}
Let $U, V, M \in \mod_\Lambda (R)$.
\begin{enumerate}

\item
If $h \in \Hom_\Lambda (U, M)$, then the homomorphism
\[
\Ext_\Lambda^1 (V, U) \to \Ext_\Lambda^1 (V, M), \; [\xi] \mapsto
[h \circ \xi],
\]
is given by
\[
Z + \bbB^{V, U} \mapsto (h_{t \alpha} Z_\alpha) + \bbB^{V, M}.
\]

\item
If $h \in \Hom_\Lambda (M, V)$, then the homomorphism
\[
\Ext_\Lambda^1 (V, U) \to \Ext_\Lambda^1 (M, U), \; [\xi] \mapsto
[\xi \circ h],
\]
is given by
\[
Z + \bbB^{V, U} \mapsto (Z_\alpha h_{s \alpha}) + \bbB^{M, U}.
\]

\end{enumerate}
\end{proposition}

\begin{proof}
We only prove the first assertion. The proof of the second one is
dual. Fix $Z \in \bbZ^{V, U}$ and let $Z_\alpha' = h_{t_\alpha}
Z_\alpha$ for $\alpha \in \Delta_1$. It follows easily that
$Z_\rho'^{V, M} = h_{t_\rho} Z_\rho^{V, U} = 0$ for all $\rho \in
\frakR$, hence $Z' \in \bbZ^{V, M}$. We define $h' \in
\Hom_\Lambda (W^Z, W^{Z'})$ by
\[
h_x' (u) = h_x (u), \; u \in U_x, \qquad \text{and} \qquad h_x'
(v) = v, \; v \in V_x,
\]
for $x \in \Delta_0$. Since $f^{Z'} h = h' f^Z$ and $g^{Z'} h' =
g^Z$ (recall that $W_x^Z = U_x \oplus V_x$ and $W_x^{Z'} = M_x
\oplus V_x$ for all $x \in \Delta_0$, $f^Z$ and $f^{Z'}$ are the
canonical injections, while $g^Z$ and $g^{Z'}$ are the canonical
projections), $[\xi^{Z'}] = [h \circ \xi^Z]$, and this finishes
the proof.
\end{proof}

We would like to have an analogous interpretation of second
extension groups and the homomorphisms
\begin{gather*}
\Ext_\Lambda^1 (V, U) \to \Ext_\Lambda^2 (V, M), \; [\xi] \mapsto
[\xi_1 \circ \xi],
\\
\intertext{and} %
\Ext_\Lambda^1 (V, U) \to \Ext_\Lambda^2 (M, U), \; [\xi] \mapsto
[\xi \circ \xi_2],
\end{gather*}
for $[\xi_1] \in \Ext_\Lambda^1 (U, M)$, $[\xi_2] \in
\Ext_\Lambda^1 (M, V)$, and $\Lambda$-modules $U$, $V$, $M$.

We first present an interpretation of second extension groups. For
$N \in \mod_\Lambda (R)$ we define $R$-representations $P^N$ and
$\Omega^N$ of $(\Delta, \frakR)$ and homomorphisms $f^N : \Omega^N
\to P^N$ and $g^N : P^N \to N$ such that $P^N$ is projective in
$\mod_\Lambda (R)$ and the sequence
\[
\xi^N : 0 \to \Omega^N \xrightarrow{f^N} P^N \xrightarrow{g^N} N
\to 0
\]
is exact, in the following way:
\begin{align*}
P_x^N & = \bigoplus_{y \in \Delta_0} x \Lambda y \otimes N_y, \; x
\in \Delta_0,
\\
\Omega_x^N & = \bigoplus_{y \in \Delta_0} x (\rad \Lambda) y
\otimes N_y, \; x \in \Delta_0,
\\
P_\alpha^N (\sigma \otimes n) & = \alpha \sigma \otimes n, \;
\alpha \in \Delta_1, \; \sigma \in s_\alpha \Lambda y, \; n \in
N_y, \; y \in \Delta_0,
\\
\Omega_\alpha^N (\sigma \otimes n) & = \alpha \sigma \otimes n -
\alpha \otimes N_\sigma n, \;
\\ 
& \qquad 
\alpha \in \Delta_1, \; \sigma \in s_\alpha (\rad \Lambda) y, \; n
\in N_y, \; y \in \Delta_0,
\\
f_x^N (\sigma \otimes n) & = \sigma \otimes n - x \otimes N_\sigma
n, \; \sigma \in x (\rad \Lambda) y, \; n \in N_y, \; x, y \in
\Delta_0,
\\
g_x^N (\sigma \otimes n) & = N_\sigma n, \; \sigma \in x \Lambda
y, \; n \in N_y, \; x, y \in \Delta_0.
\end{align*}

We now assume until the end of the section that $R = k$, i.e.\ we
are working in the category $\mod_\Lambda$.

\begin{proposition}
If $M$ and $N$ are $\Lambda$-modules, then the map
\[
\Ext_\Lambda^1 (\Omega^N, M) \to \Ext_\Lambda^2 (N, M), \; [\xi]
\mapsto [\xi \circ \xi^N],
\]
is an isomorphism.
\end{proposition}

\begin{proof}
Since $P^N$ is projective, the claim follows by standard
homological algebra (see for
example~\cite{MacLane1995}*{Chapter~III}).
\end{proof}

From now on we identify $\Ext_\Lambda^2 (N, M)$ with
$\Ext_\Lambda^1 (\Omega^N, M)$, hence with $\bbZ^{\Omega^N, M} /
\bbB^{\Omega^N, M}$, for any $\Lambda$-modules $M$ and $N$.

\begin{proposition} \label{prop_extext}
Let $U$, $V$ and $M$ be $\Lambda$-modules.
\begin{enumerate}

\item
If $Z \in \bbZ^{U, M}$, then the homomorphism
\[
\Ext_\Lambda^1 (V, U) \to \Ext_\Lambda^2 (V, M), \; [\xi] \mapsto
[\xi^Z \circ \xi],
\]
is given by $Z' + \bbB^{V, U} \mapsto Z'' + \bbB^{\Omega^V, M}$,
where
\[
Z_\alpha'' (\sigma \otimes v) = Z_\alpha Z_\sigma'^{V, U} v, \;
\alpha \in \Delta_1, \; \sigma \in s_\alpha (\rad \Lambda) y, \; v
\in V_y, \; y \in \Delta_0.
\]

\item
If $Z \in \bbZ^{M, V}$, then the homomorphism
\[
\Ext_\Lambda^1 (V, U) \to \Ext_\Lambda^2 (M, U), \; [\xi] \mapsto
[\xi \circ \xi^Z],
\]
is given by $Z' + \bbB^{V, U} \mapsto Z'' + \bbB^{\Omega^M, U}$,
where
\[
Z_\alpha'' (\sigma \otimes m) = Z_\alpha' Z_\sigma^{M, V} m, \;
\alpha \in \Delta_1, \; \sigma \in s_\alpha (\rad \Lambda) y, \; m
\in M_y, \; y \in \Delta_0.
\]

\end{enumerate}
\end{proposition}

\begin{proof}
We only prove the first assertion, the prove of the second one is
analogous. Fix $Z' \in \bbZ^{V, U}$ and let $Z'' \in
\bbA^{\Omega^V, M}$ be defined in the way described in the
proposition. One easily checks that $Z'' \in Z^{\Omega^V, M}$. We
need to show that $[\xi^Z \circ \xi^{Z'}] = [\xi^{Z''} \circ
\xi^V]$. In order to do this, it is enough to construct
homomorphisms $h : W^{Z''} \to W^Z$ and $h' : P^V \to W^{Z'}$ such
that $f^Z = h f^{Z''}$, $f^{Z'} g^Z h = h' f^V g^{Z''}$, and
$g^{Z'} h = g^V$. Such maps are defined by
\begin{align*}
h_x (m) & = m, \; m \in M_x,
\\
h_x (\sigma \otimes v) & = Z_\sigma'^{V, U} v, \; \sigma \in x
(\rad \Lambda) y, \; v \in V_y, \; y \in \Delta_0,
\\
\intertext{and} %
h_x' (\sigma \otimes v) & = Z_\sigma'^{V, U} v + V_\sigma v, \;
\sigma \in x \Lambda y, \; v \in V_y, \; y \in \Delta_0,
\end{align*}
for $x \in \Delta_0$.
\end{proof}

We now construct a ``smaller'' model of $\Ext_\Lambda^2 (N, M)$
for $M, N \in \mod_\Lambda$. Let $\bbB'^{N, M}$ be the image of
the map $\bbA^{N, M} \to \bbR^{N, M}$, $Z \mapsto (Z_\rho^{N,
M})$. We define also $\Phi^{N, M} : \bbZ^{\Omega^N, M} \to
\bbR^{N, M}$ by
\[
\Phi^{N, M} (Z)_{\rho} (n) = \sum_{i \in [1, l]} \sum_{j \in [1,
m_i - 1]} \lambda_i M_{\alpha_{i, 1} \cdots \alpha_{i, j - 1}}
Z_{\alpha_{i, j}} (\alpha_{i, j + 1} \cdots \alpha_{i, m_i}
\otimes n)
\]
for $\rho \in \frakR$ and $n \in N_{s_\rho}$, if $\rho = \lambda_1
\sigma_1 + \cdots + \lambda_l \sigma_l$ for scalars $\lambda_1,
\ldots, \lambda_l \in k$ and paths $\sigma_i = \alpha_{i, 1}
\cdots \alpha_{i, m_i} \in t_\rho (k \Delta) s_\rho$ with
$\alpha_{i, 1}, \ldots, \alpha_{i, m_i} \in \Delta_1$, $i \in [1,
l]$.

\begin{proposition}
Assume that $\Lambda$ is tilted. If $M, N \in \mod_\Lambda$, then
$\Phi^{N, M}$ is an epimorphism such that
\[
\Ker \Phi^{N, M} \subset \bbB^{\Omega^N, M} \qquad \text{and}
\qquad \Phi^{N, M} (\bbB^{\Omega^N, M}) = \bbB'^{N, M}.
\]
In particular, $\Phi^{N, M}$ induces an isomorphism
\[
\Ext_\Lambda^2 (N, M) \simeq \bbR^{N, M} / \bbB'^{N, M}.
\]
\end{proposition}

\begin{proof}
We first show that $\Ker \Phi^{N, M} \subset \bbB^{\Omega^N, M}$,
i.e.\ for $Z \in \Ker \Phi^{N, M}$ we construct $h \in
\bbV^{\Omega^N, M}$ such that $Z_\alpha = M_\alpha h_{s_\alpha} -
h_{t_\alpha} \Omega_\alpha^N$ for all $\alpha \in \Delta_1$. The
definition of such $h$ is recursive, namely
\begin{align*}
h_{t_\alpha} (\alpha \otimes n) & = 0, \; \alpha \in \Delta_1, \;
n \in N_{s_\alpha},
\\
\intertext{and} %
h_{t_\alpha} (\alpha \sigma \otimes n) & = M_\alpha h_{s_\alpha}
(\sigma \otimes n) - Z_\alpha (\sigma \otimes n), \;
\\ 
& \qquad 
\alpha \in \Delta_1, \; \sigma \in s_\alpha (\rad \Lambda) y, \; n
\in N_y, \; y \in \Delta_0.
\end{align*}

In order to prove that $\Phi^{N, M} (\bbB^{\Omega^N, M}) \subset
\bbB'^{N, M}$, take $h \in \bbV^{\Omega^N, M}$ and let $Z_\alpha =
M_\alpha h_{s_\alpha} - h_{t_\alpha} \Omega_\alpha^N$ for $\alpha
\in \Delta_1$. Direct calculations show that
\[
\Phi^{N, M} (Z)_\rho (n) = \sum_{i \in [1, l]} \sum_{j \in [1, m_i
- 1]} \lambda_i M_{\alpha_{i, 1} \cdots \alpha_{i, j - 1}}
h_{t_{\alpha_{i, j}}} (\alpha_{i, j} \otimes N_{\alpha_{i, j + 1}
\cdots \alpha_{i, m_i}} n)
\]
for $\rho$ as above (in the definition of $\Phi^{N, M}$) and $n
\in N_{s_\rho}$, i.e.\ $\Phi^{N, M} (Z) = (Z_\rho'^{N, M})$ for
$Z' \in \bbA^{N, M}$ defined by $Z'_\alpha (n) = h_{t_\alpha}
(\alpha \otimes n)$ for $\alpha \in \Delta_1$ and $n \in N_{s
\alpha}$.

We now show that $\bbB'^{N, M} \subset \Phi^{N, M}
(\bbB^{\Omega^N, M})$. Take $Z \in \bbA^{N, M}$. The above
calculations show that if we define $h \in \bbV^{\Omega^N, M}$ by
\begin{align*}
h_{t_\alpha} (\alpha \otimes n) & = Z_\alpha (n), \; \alpha \in
\Delta_1, \; n \in N_{s_\alpha},
\\
\intertext{and} %
h_{t_\alpha} (\alpha \sigma \otimes n) & = 0, \; \alpha \in
\Delta_1, \; \sigma \in s_\alpha (\rad \Lambda) y, \; n \in N_y,
\; y \in \Delta_0,
\end{align*}
and $Z' \in \bbB^{\Omega^N, M}$ by $Z'_\alpha = M_\alpha
h_{s_\alpha} - h_{t_\alpha} \Omega_\alpha^N$ for $\alpha \in
\Delta_1$, then $\Phi^{N, M} (Z') = (Z_\rho^{N, M})$.

Finally we show that $\Im \Phi^{N, M} = \bbR^{N, M}$. Observe that
$\dim_k \bbR^{N, M} = \sum_{\rho \in R} d_{s_\rho}'' d_{t_\rho}'$,
where $\bd' = \bdim M$ and $\bd'' = \bdim N$. On the other hand,
\begin{align*}
\lefteqn{\dim_k \Im \Phi^{N, M}} \quad & \\
& = \dim_k \bbZ^{\Omega^N, M} - \dim_k \Ker \Phi^{N, M}
\\
& = (\dim_k \bbZ^{\Omega^N, M} - \dim_k \bbB^{\Omega^N, M}) +
(\dim_k \bbB^{\Omega^N, M} - \dim_k \Ker \Phi^{N, M})
\\
& = [N, M]^2 + \dim_k \bbB'^{N, M}
\\
& = [N, M]^2 + \dim_k \bbA^{N, M} - \dim_k \bbZ^{N, M}
\\
& = [N, M]^2 + \sum_{\alpha \in \Delta_1} d_{s_\alpha}''
d_{t_\alpha}' - [N, M]^1 + [N, M] - \sum_{x \in \Delta_0} d_x''
d_x'
\\
& = \langle \bd'', \bd' \rangle + \sum_{\alpha \in \Delta_1}
d_{s_\alpha}'' d_{t_\alpha}' - \sum_{x \in \Delta_0} d_x'' d_x' =
\sum_{\rho \in R} d_{s_\rho}'' d_{t_\rho}',
\end{align*}
hence the claim follows.
\end{proof}

We remark that we could replace the assumption $\Lambda$ is tilted
by a more general assumption that there are no oriented cycles in
$\Delta$ and $\gldim \Lambda \leq 2$. Moreover, this assumption is
necessary only for proving that $\Phi^{N, M}$ is surjective, or
more generally, in order to identify $\Im \Phi^{N, M}$. We also
mention, that the above obtained description of second extension
groups is a trace of a much more general construction investigated
by Butler and King~\cite{ButlerKing1999}. In particular, using
their results one could identify $\Im \Phi^{N, M}$ in a general
case.

Using the above isomorphism we can give another formulation of
Proposition~\ref{prop_extext}. This formulation could also be
deduced from~\cite{ButlerKing1999}. Let $U$, $V$ and $W$ be
$\Lambda$-modules. For $Z' \in \bbZ^{V, U}$ and $Z'' \in \bbZ^{W,
V}$ we define $Z' \circ Z'' \in \bbR^{W, U}$ by
\begin{multline*} 
(Z' \circ Z'')_\rho =
\\ 
\sum_{i \in [1, l]} \sum_{\substack{j_1, j_2 \in [1, m_i - 1] \\
j_1 < j_2}} \lambda_i U_{\alpha_{i, 1} \cdots \alpha_{i, j_1 - 1}}
Z'_{\alpha_{i, j_1}} V_{\alpha_{i, j_1 + 1} \cdots \alpha_{i, j_2
- 1}} Z''_{\alpha_{i, j_2}} W_{\alpha_{i, j_2 + 1} \cdots
\alpha_{i, m_i}}
\end{multline*} 
for $\rho \in \frakR$, if $\rho = \lambda_1 \sigma_1 + \cdots +
\lambda_l \sigma_l$ with scalars $\lambda_1, \ldots, \lambda_l \in
k$ and paths $\sigma_i = \alpha_{i, 1} \cdots \alpha_{i, m_i} \in
t_\rho (k \Delta) s_\rho$, where $\alpha_{i, 1}, \ldots,
\alpha_{i, m_i} \in \Delta_1$, $i \in [1, l]$.

\begin{proposition} \label{prop_extextbis}
Let $U$, $V$ and $M$ be $\Lambda$-modules.
\begin{enumerate}

\item
If $Z \in \bbZ^{U, M}$, then, under the isomorphism induced by
$\Phi^{V, M}$, the homomorphism
\[
\Ext_\Lambda^1 (V, U) \to \Ext_\Lambda^2 (V, M), \; [\xi] \mapsto
[\xi^Z \circ \xi],
\]
is given by $Z' + \bbB^{V, U} \mapsto Z \circ Z' + \bbB'^{V, M}$.

\item
If $Z \in \bbZ^{M, V}$, then, under the isomorphism induced by
$\Phi^{M, U}$, the homomorphism
\[
\Ext_\Lambda^1 (V, U) \to \Ext_\Lambda^2 (M, U), \; [\xi] \mapsto
[\xi \circ \xi^Z],
\]
is given by $Z' + \bbB^{V, U} \mapsto Z' \circ Z + \bbB'^{M, U}$.

\end{enumerate}
\end{proposition}

\begin{proof}
Direct calculations.
\end{proof}

\section{Module schemes} \label{sect_scheme}

Throughout this section $\Lambda$ is the path algebra of a bound
quiver $(\Delta, \frakR)$. We define in this section the schemes
of $\Lambda$-modules. We also investigate subschemes of products
of module schemes consisting of pairs of modules with given
dimensions of the homomorphism space and the first extension
group.

For $\bd \in \bbN^{\Delta_0}$ and a commutative $k$-algebra $R$,
let $R^{\bd} = (R^{d_x})_{x \in \Delta_0}$. If $M \in
\bbA^{R^{\bd}, R^{\bd}}$, then we can treat $M$ as an
$R$-representation of $\Delta$ by taking $M_x = R^{d_x}$ for $x
\in \Delta_0$. In particular, $M_\rho$ is defined for each $\rho
\in x (k \Delta) y$ with $x, y \in \Delta_0$. Obviously, if $M \in
\bbA^{R^{\bd}, R^{\bd}}$, then $M \in \mod_\Lambda (R)$ if and
only if $M_\rho = 0$ for all $\rho \in R$. We define
$\mod_\Lambda^{\bd} (R)$ as the subset of all elements of
$\bbA^{R^{\bd}, R^{\bd}}$ for which the above condition is
satisfied. Then $\mod_\Lambda^{\bd}$ is a functor from the
category of commutative $k$-algebras to the category of sets,
which is an affine scheme. We call $\mod_\Lambda^{\bd}$ the scheme
of $\Lambda$-modules of dimension vector $\bd$. Note that if $M
\in \mod_\Lambda (R)$ and $\bdim M = \bd$, then there exists
(usually nonunique) $M' \in \mod_\Lambda^{\bd} (R)$ such that $M
\simeq M'$. Consequently, we will usually treat $\Lambda$-modules
as elements of $\mod_\Lambda^{\bd} (k)$.

We now generalize a construction described by
Zwara~\cite{Zwara2002a}*{Section~3}. For $\bd', \bd'' \in
\bbN^{\Delta_0}$ and $d \in \bbN$ we denote by $\calH^{\bd',
\bd''}_d$ the subscheme of $\mod_\Lambda^{\bd'} \times
\mod_\Lambda^{\bd''}$ such that $\calH^{\bd', \bd''}_d (k)$
consists of all $(U, V) \in \mod_\Lambda^{\bd'} (k) \times
\mod_\Lambda^{\bd''} (k)$ such that $[V, U] = d$. We briefly
describe its construction.

Fix a commutative $k$-algebra $R$. Recall that for $V \in
\mod_\Lambda^{\bd''} (R)$ we constructed in Section~\ref{sect_ext}
the exact sequence
\[
0 \to \Omega^V \xrightarrow{f^V} P^V \xrightarrow{g^V} V \to 0
\]
with $P^V$ projective. Iterating this construction we obtain the
exact sequence
\[
P^{\Omega^V} \xrightarrow{f^V g^{\Omega^V}} P^V \xrightarrow{g^V}
V \to 0.
\]
In particular, for each $U \in \mod_\Lambda^{\bd'} (R)$ we have
the exact sequence
\begin{scriptsize} 
\[
0 \to \Hom_\Lambda (V, U) \xrightarrow{\Hom_\Lambda (g^V, U)}
\Hom_\Lambda (P^V, U) \xrightarrow{\Hom_\Lambda (f^V g^{\Omega^V},
U)} \Hom_\Lambda (P^{\Omega^V}, U).
\]
\end{scriptsize} 
Recall that
\[
P^V = \bigoplus_{x \in \Delta_0} \Lambda x \otimes R^{d_x''}
\qquad \text{and} \qquad P^{\Omega^V} = \bigoplus_{x, y \in
\Delta_0} \Lambda x \otimes x (\rad \Lambda) y \otimes R^{d_y''}.
\]
Consequently,
\begin{align*}
\Hom_\Lambda (P^V, U) & \simeq \bigoplus_{x \in \Delta_0} R^{d_x'
d_x''}
\\ 
\intertext{and} 
\Hom_\Lambda (P^{\Omega^V}, U) & \simeq \bigoplus_{x, y \in
\Delta_0} R^{d_x' d_y''} \otimes x (\rad \Lambda) y.
\end{align*}
Moreover,
\[
f_z^V g_z^{\Omega^V} (\sigma_1 \otimes \sigma_2 \otimes v) =
\sigma_1 \sigma_2 \otimes v - \sigma_1 \otimes V_{\sigma_2} v
\]
for all $\sigma_1 \in z \Lambda x$, $\sigma_2 \in x (\rad \Lambda)
y$, $v \in R^{d_y''}$, $x, y, z \in \Delta_0$, hence if we fix
bases in $x (\rad \Lambda) y$ for all $x, y \in \Delta_0$, then
simple (although technical) analysis shows that there exists a
morphism $\varphi : \mod_\Lambda^{\bd'} \times
\mod_\Lambda^{\bd''} \to \bbM_{p, q}$ such that $\Hom_\Lambda (f^V
g^{\Omega^V}, U)$ is given by $\varphi (U, V)$ for all $U \in
\mod_\Lambda^{\bd'}$ and $V \in \mod_\Lambda^{\bd''}$, where for
\[
p = \sum_{x \in \Delta_0} d_x' d_x'' \qquad \text{and} \qquad q =
\sum_{x, y \in \Delta_0} d_x' d_y'' \dim_k (x (\rad \Lambda) y)
\]
$\bbM_{p, q}$ denotes the scheme of $p \times q$-matrices. We
define $\calH_d^{\bd', \bd''}$ to be the inverse image by
$\varphi$ of the reduced subscheme of $\bbM_{p, q}$ whose
$k$-rational points are the $p \times q$-matrices of rank $r$ with
coefficients in $k$ (see~\cite{Zwara2002a}*{Subsection~3.2}).
Obviously,
\[
\calH_d^{\bd', \bd''} (k) = \{ (U, V) \in \mod_\Lambda^{\bd'} (k)
\times \mod_\Lambda^{\bd''} (k) \mid [V, U] = d \}.
\]

Let $R_0 = k [\varepsilon] / (\varepsilon^2)$ and let $\pi : R_0
\to k$ be the canonical projection. Recall, that if $\calF$ is a
scheme, then for $x \in \calF (k)$, $\calF (\pi)^{-1} (x)$ can be
interpreted as the tangent space $T_x \calF (k)$ to $\calF (k)$ at
$x$. The proof of the following proposition just repeats the
arguments used in the proof of~\cite{Zwara2002a}*{Lemma~3.5},
hence we omit it.

\begin{lemma}
Let $M \in \mod_\Lambda^{\bd'} (R_0)$ and $N \in
\mod_\Lambda^{\bd''} (R_0)$. Put $U = \mod_\Lambda^{\bd'} (\pi)
(M)$ and $V = \mod_\Lambda^{\bd''} (\pi) (N)$. Then $(M, N) \in
\calH_d^{\bd', \bd''} (R_0)$ if and only if $[V, U] = d$ and $[N,
M] = 2 d$.
\end{lemma}

We now give another interpretation of this result. Fix $W \in
\mod_\Lambda^{\bd} (k)$ for a dimension vector $\bd$. If $L \in
\mod_\Lambda^{\bd} (\pi)^{-1} (W)$, then $L = W + \varepsilon
\ol{L}$ for some $\ol{L} \in \bbA^{k^{\bd}, k^{\bd}}$. Observe
that
\[
L_\rho = W_\rho + \varepsilon \ol{L}_\rho^{W, W}
\]
for all $\rho \in x (k \Delta) y$ with $x, y \in \Delta_0$. This
implies that the map
\[
\mod_\Lambda^{\bd} (\pi)^{-1} (W) \ni L \mapsto \ol{L} \in
\bbZ^{W, W}
\]
is well-defined and bijective.

\begin{proposition}
Let $M \in \mod_\Lambda^{\bd'} (R_0)$ and $N \in
\mod_\Lambda^{\bd''} (R_0)$. Put $U = \mod_\Lambda^{\bd'} (\pi)
(M)$ and $V = \mod_\Lambda^{\bd''} (\pi) (N)$. Then $(M, N) \in
\calH_d^{\bd', \bd''} (R_0)$ if and only if $[V, U] = d$ and
$[\xi^{\ol{M}} \circ f] = [f \circ \xi^{\ol{N}}]$ for all $f \in
\Hom_\Lambda (V, U)$.
\end{proposition}

\begin{proof}
We need to show that under the assumption $[V, U] = d$, the
condition $[N, M] = 2 d$ is equivalent to the latter condition in
the proposition. Observe that a homomorphism $f : N \to M$ is of
the form $f = f^0 + \varepsilon f^1$ for $f^0, f^1 \in
\bbV^{k^{\bd''}, k^{\bd'}}$ such that
\[
U_\alpha f_{s_\alpha}^0 = f_{t_\alpha}^0 V_\alpha \qquad
\text{and} \qquad U_\alpha f_{s_\alpha}^1 + \ol{M}_\alpha
f_{s_\alpha}^0 = f_{t_\alpha}^1 V_\alpha + f_{t_\alpha}^0
\ol{N}_\alpha
\]
for all $\alpha \in \Delta_1$. The first condition means that $f^0
\in \Hom_\Lambda (V, U)$, hence the dimension of the set of such
$f^0$ equals $d$. Consequently, the condition $[N, M] = 2 d$ is
equivalent to the statement that for each $f \in \Hom_\Lambda (V,
U)$ the set of $g \in \bbV^{k^{\bd'}, k^{\bd''}}$ such that
\[
\ol{M}_\alpha f_{s_\alpha} - f_{t_\alpha} \ol{N}_\alpha =
g_{t_\alpha} V_\alpha - U_\alpha g_{s_\alpha}
\]
for all $\alpha \in \Delta_1$, has dimension $d$. However linear
algebra says that this is that case if and only if this system of
equations has a solution, i.e.\ if and only if $(\ol{M}_\alpha
f_{s_\alpha} - f_{t_\alpha} \ol{N}_\alpha) \in \bbB^{V, U}$, what
finishes the proof according to Proposition~\ref{prop_exthom}.
\end{proof}

We now extend the above construction to extensions. Fix $d' \in
\bbN$. If $V \in \mod_\Lambda^{\bd''} (k)$, then we have the exact
sequence
\[
P^{\Omega^{\Omega^V}} \xrightarrow{f^{\Omega^V}
g^{\Omega^{\Omega^V}}} P^{\Omega^V} \xrightarrow{f^V g^{\Omega^V}}
P^V \xrightarrow{g^V} V \to 0.
\]
Repeating the previous construction (using also $\Hom_\Lambda
(f^{\Omega^V} g^{\Omega^{\Omega^V}}, U)$ this time) we define the
subscheme $\calE_{d, d'}^{\bd', \bd''}$ of $\calH_d^{\bd', \bd''}$
such that $(U, V) \in \calE_{d, d'}^{\bd', \bd''} (k)$ for $(U, V)
\in \calH_d^{\bd', \bd''} (k)$ if and only if $[V, U]^1 = d'$.
Moreover, if $(M, N) \in \calH_d^{\bd', \bd''} (R_0)$, $U =
\mod_\Lambda^{\bd'} (\pi) (M)$ and $V = \mod_\Lambda^{\bd''} (\pi)
(N)$, then $(M, N) \in \calE_{d, d'}^{\bd', \bd''} (R_0)$ if and
only if $[V, U]^1 = d'$ and $[N, M]^1 = 2 d'$. An alternate
description of $\calE_{d, d'}^{\bd', \bd''} (R_0)$ is the
following.

\begin{proposition} \label{prop_schemeext}
Let $M \in \mod_\Lambda^{\bd'} (R_0)$ and $N \in
\mod_\Lambda^{\bd''} (R_0)$. Put $U = \mod_\Lambda^{\bd'} (\pi)
(M)$ and $V = \mod_\Lambda^{\bd''} (\pi) (N)$. Then $(M, N) \in
\calE_{d, d'}^{\bd', \bd''} (R_0)$ if and only if $(M, N) \in
\calH_d^{\bd', \bd''} (R_0)$, $[V, U]^1 = d'$, and $[\xi^{\ol{M}}
\circ \xi] + [\xi \circ \xi^{\ol{N}}] = 0$ for all $[\xi] \in
\Ext_\Lambda^1 (V, U)$.
\end{proposition}

\begin{proof}
We need to show that, under the assumptions $[V, U]^1 = d'$ and
$(M, N) \in \calH_d^{\bd', \bd''} (R_0)$, the condition $[N, M]^1
= 2 d'$ is equivalent to the condition $[\xi^{\ol{M}} \circ \xi] +
[\xi \circ \xi^{\ol{N}}] = 0$ for all $[\xi] \in \Ext_\Lambda^1
(V, U)$. Thus assume that $(M, N) \in \calH_d^{\bd', \bd''} (R_0)$
and $[V, U] = d'$. Since $\Ext_\Lambda^1 (V, U) = \bbZ^{V, U} /
\bbB^{V, U}$ and $\Ext_\Lambda^1 (N, M) = \bbZ^{N, M} / \bbB^{N,
M}$, it follows that $[N, M]^1 = 2 d' = 2 [V, U]^1$ if and only if
$\dim_k \bbZ^{N, M} = 2 \dim_k \bbZ^{V, U}$.

If $Z \in \bbA^{N, M}$, then $Z = Z' + \varepsilon \ol{Z}$ for
some $Z', \ol{Z} \in \bbA^{k^{\bd''}, k^{\bd'}}$. Moreover, easy
calculations show that $Z \in \bbZ^{N, M}$, i.e.\ $Z_\rho^{N, M} =
0$ for all $\rho \in \frakR$, if and only if $Z' \in \bbZ^{V, U}$
and
\[
\ol{M} \circ Z' + Z' \circ \ol{N} + \Phi^{V, U} (\ol{Z}) = 0.
\]
Consequently, $\dim_k \bbZ^{N, M} = 2 \dim_k \bbZ^{V, U}$ if and
only if $Z' \in \bbZ^{V, U}$ and the set of $\ol{Z} \in
\bbA^{k^{\bd''}, k^{\bd'}}$ such that the above condition is
satisfied has dimension $\dim_k \bbZ^{V, U}$. Since $\Ker \Phi^{V,
U} = \bbZ^{V, U}$, linear algebra says that $\dim_k \bbZ^{N, M} =
2 \dim_k \bbZ^{V, U}$ if and only if $\ol{M} \circ Z' + Z' \circ
\ol{N} \in \Im \Phi^{V, U} = \bbB'^{V, U}$ for all $Z' \in
\bbZ^{V, U}$, and this finishes the proof according to
Proposition~\ref{prop_extextbis}.
\end{proof}

Let $\bd \in \bbN^{\Delta_0}$. The product $\GL_{\bd} (k) =
\prod_{x \in \Delta_0} \GL_{d_x} (k)$ of general linear groups
acts on $\mod_\Lambda^{\bd}$ by conjugation:
\[
(g \cdot M)_\alpha = g_{t_\alpha} M_\alpha g_{s_\alpha}^{-1}, \; g
\in \GL_{\bd} (k), \; M \in \mod_\Lambda^{\bd} (k), \; \alpha \in
\Delta_1.
\]
If $\calO (M)$ denotes the $\GL_{\bd} (k)$-orbit of $M \in
\mod_\Lambda^{\bd} (k)$, then $\calO (M_1) = \calO (M_2)$ for
$M_1, M_2 \in \mod_\Lambda^{\bd} (k)$ if and only if $M_1 \simeq
M_2$. It is known (see for
example~\cite{KraftRiedtmann1986}*{2.2}) that
\[
\dim \calO (M) = \dim \GL_{\bd} (k) - [M, M]
\]
for $M \in \mod_\Lambda^{\bd} (k)$. We put
\[
a (\bd) = \dim \bbA^{k^{\bd}, k^{\bd}} - \sum_{\rho \in \frakR}
d_{s_\rho} d_{t_\rho}.
\]
Note that $a (\bd) = \dim \GL_{\bd} (k) - \chi (\bd)$. Fix $M \in
\mod_\Lambda^{\bd} (k)$. It follows from Voigt's
result~\cite{Voigt1977} that the Zariski closure $\ol{\calO (M)}$
of $\calO (M)$ is an irreducible component of $\mod_\Lambda^{\bd}
(k)$ provided $[M, M]^1 = 0$. Moreover, if $[M, M]^2 = 0$, then
$M$ is a regular point of $\mod_\Lambda^{\bd}
(k)$~\cite{Geiss1996}, i.e.\ $\dim_k T_M \mod_\Lambda^{\bd} (k)$
equals the maximum of the dimensions of the irreducible components
of $\mod_\Lambda^{\bd} (k)$ containing $M$. Finally, if $\Lambda$
is tilted and $[M, M]^1 = 0 = [M, M]^2$, then $\dim \ol{\calO (M)}
= a (\bd)$.

Let $M, N \in \mod_\Lambda^{\bd} (k)$ for $\bd \in
\bbN^{\Delta_0}$. We call $N$ a degeneration of $M$ and write $M
\leq_{\deg} N$ if $N\in \ol{\calO (M)}$. A deneneration $M
\leq_{\deg} N$ is called minimal if $M \not \simeq N$ and either
$M \simeq L$ or $N \simeq L$ for each $L \in \mod_\Lambda^{\bd}
(k)$ such that $M \leq_{\deg} L \leq_{\deg} N$. If $M \leq_{\deg}
N$ and $M \not \simeq N$, then we write $M <_{\deg} N$.

\section{Proof of the main result} \label{sect_proof}

Throughout this section we fix a sincere directing module $M$ over
the path algebra $\Lambda$ of a bound quiver $(\Delta, \frakR)$.
Recall, that this implies that $\Lambda$ is a tilted algebra. We
also fix a directing tilting module $T$ such that $\add M \subset
\add T$. Let $T_1$, \ldots, $T_n$ be the pairwise nonisomorphic
indecomposable direct summands of $T$, $\calF = \calF (T)$ and
$\calT = \calT (T)$. Upper semicontinuity of the functions $[T,
-], [-, T]^1 : \mod_\Lambda^{\bd} (k) \to \bbZ$ (see for
example~\cite{Crawley-BoeveySchroer2002}*{Lemma~4.3}) implies that
$\calF^{\bd}$ and $\calT^{\bd}$ are open subset of
$\mod_\Lambda^{\bd} (k)$ for each $\bd \in \bbN^{\Delta_0}$, where
$\calA^{\bd} = \calA \cap \mod_\Lambda^{\bd} (k)$ for a
subcategory $\calA$ of $\mod_\Lambda$ and $\bd \in
\bbN^{\Delta_0}$. In order to generalize the above observation we
need the following lemma.

\begin{lemma}
If $N', N'' \in \calT$ and $[T_i, N'] = [T_i, N'']$ for all $i \in
[1, n]$, then $\bdim N' = \bdim N''$.
\end{lemma}

\begin{proof}
We have the isomorphism $\Phi : \bbZ^{\Delta_0} \to \bbZ^n$
induced by the assignment $\bdim N \mapsto ([T_i, N] - [T_i,
N]^1)$ (see~\cite{AssemSimsonSkowronski2006}*{Theorem~VI.4.3}).
Since $\Phi (\bdim N) = ([T_i, N])$ for all $N \in \calT$, this
implies the claim.
\end{proof}

For $\calU' \subset \mod_\Lambda^{\bd'} (k)$ and $\calU'' \subset
\mod_\Lambda^{\bd''} (k)$, where $\bd', \bd'' \in
\bbN^{\Delta_0}$, let $\calU' \oplus \calU''$ be the subset of
$\mod_\Lambda^{\bd' + \bd''} (k)$ consisting of $M$ such that $M
\simeq M' \oplus M''$ for some $M' \in \calU'$ and $M'' \in
\calU''$.

\begin{corollary} \label{coro_FT}
If $\bd', \bd'' \in \bbN^{\Delta_0}$, then $\calF^{\bd'} \oplus
\calT^{\bd''}$ is a locally closed subset of $\mod_\Lambda^{\bd' +
\bd''} (k)$.
\end{corollary}

\begin{proof}
Let $\bd = \bd' + \bd''$ and $N \in \mod_\Lambda^{\bd} (k)$. If $N
\simeq N' \oplus N''$ for $N' \in \calF$ and $N'' \in \calT$, then
$[T_i, N] = [T_i, N''] = \langle \bdim T_i, \bdim N'' \rangle$ for
each $i \in [1, n]$. Consequently, it follows from the previous
lemma that $N \in \calF^{\bd'} \oplus \calT^{\bd''}$, i.e.\ $\bdim
N'' = \bd''$, if and only if $[T_i, N] = \langle \bdim T_i, \bd''
\rangle$ for all $i \in [1, n]$, hence the claim follows.
\end{proof}

\begin{lemma}
If $N$ is a minimal degeneration of $M$, then there exists an
exact sequence
\[
0 \to U \to M \to V \to 0
\]
such that $N \simeq U \oplus V$.
\end{lemma}

\begin{proof}
If there is not such a sequence, then the minimality of the
degeneration $M <_{\deg} N$ and~\cite{Zwara2000}*{Theorem~4} imply
that there exists an indecomposable direct summand $X$ of $N$ such
that $M' <_{\deg} X$ for a direct summand $M'$ of $M$. In
particular,
\[
[M', X] = [M', X]^1 + \chi (\bdim M') > 0.
\]
Similarly, $[X, M'] \neq 0$. Consequently, $X \in \add T$, hence
$X \simeq M'$, a contradiction.
\end{proof}

\begin{corollary} \label{coro_mindeg}
If $U \in \calF$ and $T \in \calT$ are such that $U \neq 0$ and $U
\oplus V$ is a minimal degeneration of $M$, then there exists an
exact sequence
\[
0 \to U \to M \to V \to 0.
\]
\end{corollary}

\begin{proof}
According to the above lemma there exists an exact sequence
\[
\sigma : 0 \to U' \to M \to V' \to 0
\]
such that $U' \oplus V' \simeq U \oplus V$. Since $[M, U] = 0$,
$U$ must be a direct summand of $U'$. Let
\[
\sigma' = p \circ \sigma : 0 \to U \to M' \to V' \to 0,
\]
where $p : U' \to U$ is the canonical projection. Since $U \not
\in \add M$, $\sigma'$ does not split. Consequently,
\[
M \leq_{\deg} \Ker p \oplus M' <_{\deg} \Ker p \oplus U \oplus V'
\simeq U' \oplus V' \simeq U \oplus V,
\]
hence $M \simeq \Ker p \oplus M'$ and $\Ker p \oplus V' \simeq V$
by the minimality of the degeneration $M <_{\deg} N$. Moreover, we
have a short exact sequence
\[
0 \to U \to M' \oplus \Ker p \to V' \oplus \Ker p \to 0
\]
and the claim follows.
\end{proof}

\begin{proof}[Proof of the main result] 
Let $\calX$ be an irreducible component of $\ol{\calO (M)}
\setminus \calO (M)$. Our aim is to show that there exists an open
subset of $\calX$ such that all points of this subset are regular
points of $\ol{\calO (M)}$. Since $\mod_\Lambda = \calF \vee
\calT$, it follows from Corollary~\ref{coro_FT} that there exist
$\bd', \bd'' \in \bbN^{\Delta_0}$ such that $\calU = (\calF^{\bd'}
\oplus \calT^{\bd''}) \cap \calX$ is an open subset of $\calX$. If
$\bd' = 0$, then $\id_\Lambda N \leq 1$ for all $N \in \calU$,
hence the claim follows. Thus assume that $\bd' \neq 0$. Let $d_0
= \min \{ [N, N] \mid N \in \calX \}$, $d_1 = \min \{ [N, N]^1
\mid N \in \calX \}$, and $\calU_0$ be the subset of $\calU$
consisting of all $N \in \calU$ such that $[N, N] = d_0$, $[N,
N]^1 = d_1$, and $N$ does not belong to an irreducible component
of $\ol{\calO (M)} \setminus \calO (M)$ different from $\calX$.
Then $\calU_0$ is an open subset of $\calX$. Moreover, all points
of $\calU_0$ are minimal degenerations of $M$. We show that all
points of $\calU_0$ are regular points of $\ol{\calO (M)}$.

We first observe that in order to prove the above claim it is
enough to show that $\dim_k L (U, V) \leq a (\bd') + a (\bd'') -
[V, U]^2$ for all $U \in \calF^{\bd'}$ and $V \in \calT^{\bd''}$
such that $U \oplus V \in \calU_0$, where
\[
L (U, V) = T_{U \oplus V} ((\mod_\Lambda^{\bd'} (k) \times
\mod_\Lambda^{\bd''} (k)) \cap \ol{\calO (M)}),
\]
for $U \in \mod_\Lambda^{\bd'} (k)$ and $V \in
\mod_\Lambda^{\bd''} (k)$, and we identify $\mod_\Lambda^{\bd'}
(k) \times \mod_\Lambda^{\bd''} (k)$ with
\[
\left\{
\begin{bmatrix}
U & 0 \\ 0 & V
\end{bmatrix}
\Big| \; U \in \mod_\Lambda^{\bd'} (k), \; V \in
\mod_\Lambda^{\bd''} (k) \right\} \subset \mod_\Lambda^{\bd} (k).
\]
Indeed, if $N \in \calU_0$, then without loss of generality we may
assume that $N = U \oplus V$ for $U \in \calF^{\bd'}$ and $V \in
\calT^{\bd''}$. Consequently,
\[
T_N \ol{\calO (M)} \subset L (U, V) \oplus \bbZ^{U, V} \oplus
\bbZ^{V, U}
\]
(here we use that $T_N \mod_\Lambda^{\bd} (k) = \bbZ^{N, N}$).
Since
\begin{align*}
\dim_k \bbZ^{U, V} & = \sum_{x \in \Delta_0} d_x' d_x'' - \langle
\bd', \bd'' \rangle
\\ 
\intertext{and} 
\dim_k \bbZ^{V, U} & = \sum_{x \in \Delta_0} d_x'' d_x' - [V, U] +
[V, U]^1,
\end{align*}
the required inequality implies $\dim_k T_n \ol{\calO (M)} \leq a
(\bd) = \dim \ol{\calO (M)}$, and this will finish the proof.

Now we fix $U \in \calF^{\bd'}$ and $V \in \calT^{\bd''}$ such
that $(U, V) \in \calU_0$. Observe that $[U, V] = \langle \bd',
\bd'' \rangle$, $[U, V]^1 = 0$, and $[V, U] = 0$. Recall also that
$[U \oplus V, U \oplus V] = d_0$ and $[U \oplus V, U \oplus V]^1 =
d_1$. Finally, $[U, U]^2 = 0$ and $[V, V]^2 = 0$. Consequently,
$[V, U]^1 = d$, where
\[
d = d_1 - d_0 + \chi (\bd') + \chi (\bd'') + \langle \bd', \bd''
\rangle,
\]
thus $\calU_0 \subset \calE_{0, d}^{\bd', \bd''}$. Observe that
$\calU_0$ is open in $(\mod_\Lambda^{\bd'} (k) \times
\mod_\Lambda^{\bd''} (k)) \cap \ol{\calO (M)}$, hence $L (U, V) =
T_{U \oplus V} \calU_0 \subset T_{(U, V)} \calE_{0, d}^{\bd',
\bd''}$. Recall that according to Proposition~\ref{prop_schemeext}
$T_{(U, V)} \calE_{0, d}^{\bd', \bd''}$ consists of $(Z', Z'') \in
\bbZ^{U, U} \times \bbZ^{V, V}$ such that $[\xi^{Z'} \circ \xi] +
[\xi \circ \xi^{Z''}] = 0$ for all $[\xi] \in \Ext_\Lambda^1 (V,
U)$. Since $U \oplus V$ is a minimal degeneration of $M$, there
exists a short exact sequence of the form
\[
\xi : 0 \to U \to M \to V \to 0
\]
by Corollary~\ref{coro_mindeg}. Moreover, $\pd_\Lambda M \leq 1$,
hence the map $[V, V]^1 \to [V, U]^2$, $[\xi'] \to [\xi \circ
\xi']$, is surjective. Consequently, the map
\[
\Psi : \bbZ^{U, U} \times \bbZ^{V, V} \to \Ext_\Lambda^2 (V, U),
\; (Z', Z'') \mapsto [\xi^{Z'} \circ \xi] + [\xi \circ \xi^{Z''}],
\]
is surjective. Finally, using $T_{(U, V)} \calE_{0, d}^{\bd',
\bd''} \subset \Ker \Psi$, $\dim_k \bbZ^{U, U} = a (\bd')$, and
$\dim_k \bbZ^{V, V} = a (\bd'')$, we finish the proof.
\end{proof}

\bibsection


\end{document}